# Stable vector bundles on algebraic surfaces

Wei-Ping Li  and  Zhenbo Qin

ABSTRACT. We prove an existence result for stable vector bundles with arbitrary rank on an algebraic surface, and determine the birational structure of certain moduli space of stable bundles on a rational ruled surface.

## 1. Introduction

Let $\mathcal{M}_L(r; c_1, c_2)$ be the moduli space of $L$-stable (in the sense of Mumford-Takemoto) rank-$r$ vector bundles with Chern classes $c_1$ and $c_2$ on an algebraic surface $X$. The nonemptiness of $\mathcal{M}_L(2; 0, c_2)$ has been studied by Taubes [22], Gieseker [9], Artamkin [1], Friedman [8], Jun Li, etc. The generic smoothness of $\mathcal{M}_L(2; c_1, c_2)$ has been proved by Donaldson [6], Friedman [8] and Zuo [23]. For an arbitrary $r$ and $c_1$, Maruyama [17] proved that for any integer $s$, there exists an integer $c_2$ with $c_2 \geq s$ such that $\mathcal{M}_L(r; c_1, c_2)$ is nonempty; however, no explicit formula for the lower bound of $c_2$ was given. Using deformation theory on torsion-free sheaves, Artamkin [1] showed that if $c_2 > (r+1) \cdot \max(1, p_g)$, then the moduli space $\mathcal{M}_L(r; 0, c_2)$ is nonempty and contains a vector bundle $V$ with $h^2(X, \text{ad}(V)) = 0$ where $\text{ad}(V)$ is the trace-free sub-vector bundle of $\mathcal{E}nd(V)$. Based on certain degeneration theory, Gieseker and J. Li [10] announced the generic smoothness of the moduli space $\mathcal{M}_L(r; c_1, c_2)$.

In the first part of this paper, we determine the nonemptiness of $\mathcal{M}_L(r; c_1, c_2)$ in the most general form, and show that at least one of the components of moduli space is generically smooth. Using an explicit construction, we show the following.

**Theorem 1.1.** *For any ample divisor $L$ on $X$, there exists a constant $\alpha$ depending only on $X$, $r$, $c_1$ and $L$ such that for any $c_2 \geq \alpha$, there exists an $L$-stable rank-$r$ bundle $V$ with Chern classes $c_1$ and $c_2$. Moreover, $h^2(X, \text{ad}(V)) = 0$.*

This is proved in section 2. Our starting point is the classical Cayley-Bacharach property. A well-known result (see p.731 in [11]) says that there exists a rank-2 bundle given by an extension of $\mathcal{O}_X(L'') \otimes I_Z$ by $\mathcal{O}_X(L')$ if and only if the 0-cycle $Z$ satisfies the Cayley-Bacharach property with respect to the complete linear system $|(L'' - L' + K_X)|$, that is, any curve in $|(L'' - L' + K_X)|$ containing all but one points in $Z$ must contain

---

the remaining point. It follows that to construct a rank-$r$ bundle $V$ as an extension of

$$\bigoplus_{i=1}^{(r-1)} [\mathcal{O}_X(L_i) \otimes I_{Z_i}]$$

by $\mathcal{O}_X(L')$, we need only to make sure that $Z_i$ satisfies the Cayley-Bacharach property with respect to $|(L_i - L' + K_X)|$ for each $i$. Now, let $L$ be an ample divisor, and normalize $c_1$ such that $-rL^2 < c_1 \cdot L \leq 0$. Let $L' = c_1 - (r-1)L$ and $L_i = L$. Our main argument is that if the length of $Z_i$ is sufficiently large and if $Z_i$ is generic in the Hilbert scheme $Hilb^{\ell(Z_i)}(X)$ for each $i$, then the vector bundle $V$ is $L$-stable and

$$h^2(X, \mathrm{ad}(V)) = 0.$$

Similar construction for stable rank-2 bundles is well-known [20].

We notice that there have been extensive studies for stable rank-2 bundles on $\mathbf{P}^2$ and on a ruled surface ([3, 14, 13, 4, 5, 8, 16, 21]), and for stable bundles with arbitrary rank on $\mathbf{P}^2$ ([15, 18, 7, 1]). In the rest of this paper, we study the structure of $\mathcal{M}_L(r; c_1, c_2)$ for a suitable ample divisor $L$ on a ruled surface $X$. In section 3, we prove that $\mathcal{M}_L(r; c_1, c_2)$ is empty if $(c_1 \cdot f)$ is not divisible by $r$, and that $\mathcal{M}_L(r; \mathbf{t}f, c_2)$ is nonempty if $-r < t \leq 0$ and $c_2 \geq 2(r-1)$; moreover, we show that the restriction of any bundle in $\mathcal{M}_L(r; \mathbf{t}f, c_2)$ to the generic fiber of the ruling $\pi$ must be trivial.

In section 4, we assume that $X$ is a rational ruled surface, and verify that a generic bundle $V$ in $\mathcal{M}_L(r; tf, c_2)$ sits in an exact sequence of the form:

$$0 \to \bigoplus_{i=1}^{r} \mathcal{O}_X(-n_i f) \to V \to \bigoplus_{i=1}^{c_2} (\tau_i)_* \mathcal{O}_{f_i}(-1) \to 0 \quad (1.2)$$

where $\{f_1, \ldots, f_{c_2}\}$ are distinct fibers with $\tau_i$ being the natural embedding $f_i \hookrightarrow X$, and the integer $n_i$ is defined inductively by (4.20). The idea is a natural generalization of those in [4, 5, 8]. Since the restriction of $V$ to the generic fiber is trivial, $\pi_* V$ is a rank-$r$ bundle on $\mathbf{P}^1$; thus, we can construct $(r-1)$ exact sequences:

$$0 \to \mathcal{O}_X(-n_i f) \to V_i^{**} \to V_{i-1} \to 0$$

where $i = r, \ldots, 2$, $V_r = V$, and $V_i$ is a torsion-free rank-$i$ sheaf. By estimating the numbers of moduli of $V_i$ and $V_i^{**}$, we conclude that for a generic $V$, the sheaves $V_2, \ldots, V_r$ are all locally free, and $V_1 = \mathcal{O}_X((c_2 - n_1)f) \otimes I_Z$ where $Z$ consists of $c_2$ points lying on distinct fibers. Then, the exact sequence (1.2) follows.

In section 5, based on (1.2), we define a rational map $\Phi$ from $\mathcal{M}_L(r; tf, c_2)$ to $\mathbf{P}^{c_2}$, and show that the fiber is unirational. We thus obtain our second main result.

**Theorem 1.3.** *Let $X$ be a rational ruled surface. Assume that the moduli space $\mathcal{M}_L(r; tf, c_2)$ is nonempty where $r \geq 2$, $-r < t \leq 0$, and $L$ satisfies the condition (3.3). Then, $\mathcal{M}_L(r; tf, c_2)$ is irreducible and unirational.*

One consequence of Theorem 1.3 is that the moduli space $\mathcal{M}_L(r; 0, c_2)$ on $\mathbf{P}^2$ which is known to be irreducible [15, 7] is unirational. In fact, we shall show that any irreducible component of a nonempty moduli space on a rational surface is unirational, and determine the irreducibility and rationality in rank-3 case. Details will appear elsewhere.

**Acknowledgment.** The authors would like to thank Jun Li and Karien O'Grady for some valuable discussions. They are very grateful to the referee for useful comments and for pointing out a mistake in the previous version. The second author also would like to thank the Institute for Advanced Study at Princeton for its hospitality and its financial support through the NSF grant DMS-9100383.

**Notations and conventions**

$X$ stands for an algebraic surface over the complex number field $\mathbf{C}$. The stability of a vector bundle is in the sense of Mumford-Takemoto. Furthermore, we make no distinction between a vector bundle and its associated locally free sheaf.

$K_X =:$ the canonical divisor of $X$;

$p_g =: h^0(X, \mathcal{O}_X(K_X))$, the geometric genus of $X$;

$\ell(Z) =:$ the length of the 0-cycle $Z$ on $X$;

$Hilb^\ell(X) =:$ the Hilbert scheme parametrizing all 0-cycles of length-$\ell$ on $X$;

$r =:$ an integer larger than one;

$\mu_L(V) =: c_1(V) \cdot L/\mathrm{rank}(V)$ where $L$ is an ample divisor on $X$ and $V$ is a torsion-free sheaf on $X$.

$\mathrm{ad}(V) =: \ker(\mathrm{Tr}: \mathcal{E}nd(V) \to \mathcal{O}_X)$. Then, $\mathcal{E}nd(V) = \mathrm{ad}(V) \oplus \mathcal{O}_X$.

$[x] =:$ the integer part of the number $x$.

When $X$ is a ruled surface, we also fix the following notations.

$\pi =:$ a ruling from $X$ to an algebraic curve $C$;

$f =:$ a fiber to the ruling $\pi$;

$\sigma =:$ a section to $\pi$ such that $\sigma^2$ is the least;

$e =: -\sigma^2$;

$r_L =: b/a$ where $L \equiv (a\sigma + bf)$ and $a \neq 0$;

$\mathbf{d}f =: \pi^*(\mathbf{d})$ where $\mathbf{d}$ is a divisor on $C$. In this case, $d$ stands for degree($\mathbf{d}$);

$\mathbf{P}_K^1 =:$ the generic fiber of the ruling $\pi$.

## 2. Existence of stable bundles on algebraic surfaces

### 2.1. The Cayley-Bacharach property

Fix divisors $L', L_1, \ldots, L_{r-1}$ and reduced 0-cycles $Z_1, \ldots, Z_{r-1}$ on the algebraic surface $X$ such that $Z_i \bigcap Z_j = \emptyset$ for $i \neq j$. Put $Z = \bigcup Z_i$ and

$$W = \bigoplus_{i=1}^{(r-1)} [\mathcal{O}_X(L_i) \otimes I_{Z_i}].$$

Let $W_i$ be the obvious quotient $W/[\mathcal{O}_X(L_i) \otimes I_{Z_i}]$. It is well known that there exists an extension $e_i$ in $\text{Ext}^1(\mathcal{O}_X(L_i) \otimes I_{Z_i}, \mathcal{O}_X(L'))$ whose corresponding exact sequence

$$0 \to \mathcal{O}_X(L') \to V_i \to \mathcal{O}_X(L_i) \otimes I_{Z_i} \to 0$$

gives a bundle $V_i$ if and only if $Z_i$ satisfies the Cayley-Bacharach property with respect to the complete linear system $|(L_i - L' + K_X)|$, i.e. if a curve $D$ in $|(L_i - L' + K_X)|$ contains all but one point of $Z_i$, then $D$ contains the remaining point. Note that

$$\text{Ext}^1(W, \mathcal{O}_X(L')) = \bigoplus_{i=1}^{(r-1)} \text{Ext}^1(\mathcal{O}_X(L_i) \otimes I_{Z_i}, \mathcal{O}_X(L')).$$

In the following, we study the existence of a bundle $V$ sitting in an extension

$$0 \to \mathcal{O}_X(L') \to V \xrightarrow{\phi} W \to 0 \qquad (2.1)$$

**Proposition 2.2.** *There exists an extension $e \in \text{Ext}^1(W, \mathcal{O}_X(L'))$ whose corresponding exact sequence (2.1) gives a bundle $V$ if and only if for each $i = 1, \ldots, (r-1)$, the 0-cycle $Z_i$ satisfies the Cayley-Bacharach property with respect to $|(L_i - L' + K_X)|$.*

*Proof.* Put $e = (e_1, \ldots, e_{r-1})$ where $e_i \in \text{Ext}^1(\mathcal{O}_X(L_i) \otimes I_{Z_i}, \mathcal{O}_X(L'))$. Let $V_i$ be the subsheaf $\phi^{-1}(\mathcal{O}_X(L_i) \otimes I_{Z_i})$ of $V$. Then, $V_i$ is given by the extension $e_i$:

$$0 \to \mathcal{O}_X(L') \to V_i \to \mathcal{O}_X(L_i) \otimes I_{Z_i} \to 0.$$

Note that $V$ is locally free outside the 0-cycle $Z$ and sits in an exact sequence

$$0 \to V_i \to V \to W_i \to 0.$$

Since $W_i$ is locally free at the points in $Z_i$, we see that $V$ is locally free at the points in $Z_i$ if and only if $V_i$ is locally free at the points in $Z_i$, that is, $Z_i$ satisfies the Cayley-Bacharach property with respect to $|(L_i - L' + K_X)|$. Hence, our result follows. ∎

**Corollary 2.3.** *If $h^0(X, \mathcal{O}_X(L_i - L' + K_X) \otimes I_{Z_i - \{x\}}) = 0$ for every $i$ and for every $x \in Z_i$, then there exists a bundle $V$ sitting in the exact sequence* (2.1).

### 2.2. Construction of a rank-$r$ bundle $V$

Let $L$ be a very ample divisor on $X$, and let $V$ be a rank-$r$ bundle. Note that

$$c_1(V \otimes \mathcal{O}_X(nL)) = c_1(V) + nrL.$$

Thus, by tensoring some line bundle to $V$, we may assume that $-rL^2 < c_1(V) \cdot L \leq 0$. Without loss of generality, from now on, we fix a divisor $c_1$ with $-rL^2 < c_1 \cdot L \leq 0$.

We start with three lemmas. In these lemmas, we prove certain properties satisfied by a generic 0-cycle in the Hilbert scheme $Hilb^\ell(X)$ when $\ell$ is sufficiently large.

**Lemma 2.4.** *Let $Z$ be a generic 0-cycle $Z$ in the Hilbert scheme $Hilb^\ell(X)$.*
  (i) *If $\ell \geq h^0(X, \mathcal{O}_X(rL - c_1 + K_X))$, then $h^0(X, \mathcal{O}_X(rL - c_1 + K_X) \otimes I_Z) = 0$;*
  (ii) *If $\ell \geq p_g$, then $h^0(X, \mathcal{O}_X(K_X) \otimes I_Z) = 0$.*

*Proof.* This is straightforward. ∎

**Lemma 2.5.** *Let $\ell \geq \max(p_g, h^0(X, \mathcal{O}_X(rL - c_1 + K_X)))$. Then, a generic 0-cycle $Z'$ in the Hilbert scheme $Hilb^{\ell+1}(X)$ satisfies the Cayley-Bacharach property with respect to $|rL - c_1 + K_X|$; moreover, $h^0(X, \mathcal{O}_X(K_X) \otimes I_{Z'}) = 0$.*

*Proof.* In view of Lemma 2.4 (ii), we need only to prove the first statement. Define an open dense subset $U_\ell$ of $Hilb^\ell(X)$ such that if $Z \in U_\ell$, then $Z$ is reduced and

$$h^0(X, \mathcal{O}_X(rL + K_X - c_1) \otimes I_Z) = 0.$$

By Lemma 2.4 (i), this can be done. Define $V_\ell$ to be the open subset of $Hilb^\ell(X)$ consisting of reduced 0-cycles. Hence $U_\ell$ is an open dense subset of $V_\ell$. Define $Z^{\ell+1}$ to be the universal family in $V_{\ell+1} \times X$:

$$Z^{\ell+1} = \{\, ([Z], x) \in V_{\ell+1} \times X \mid x \in Z\,\}.$$

Then, there is a surjective morphism $\pi : Z^{\ell+1} \to V_\ell$ given by $\pi([Z], x) = (Z - x)$. Hence, $Z^{\ell+1} - \pi^{-1}(U_\ell)$ is a proper closed subset of $Z^{\ell+1}$. Define the natural projection:

$$Z^{\ell+1} \subset V_{\ell+1} \times X \xrightarrow{\rho} V_{\ell+1}.$$

Then, $\rho$ is a flat surjection, and $\rho(Z^{\ell+1} - \pi^{-1}(U_\ell))$ is a proper closed subset of $V_{\ell+1}$. So we can choose an element $Z' \in V_{\ell+1} - \rho(Z^{\ell+1} - \pi^{-1}(U_\ell))$. Hence, $\rho^{-1}([Z']) \subset \pi^{-1}(U_\ell)$; this means that for any point $x$ in $Z'$, $Z' - x \in U_\ell$, that is, we have

$$h^0(X, \mathcal{O}_X(rL + K_X - c_1) \otimes I_{Z'-x}) = 0 \qquad \text{for any } x \in Z'.$$

So $Z'$ satisfies Cayley-Bacharach property with respect to $|rL + K_X - c_1|$. ∎

The above two lemmas will be used to construct a rank-$r$ bundle while the following lemma will be used to show the $L$-stability of that bundle.

**Lemma 2.6.** *There exists a reduced 0-cycle $Z''$ of length $\ell(Z'') \geq 4(r-1)^2 \cdot L^2$ such that if $h^0(X, \mathcal{O}_X(F) \otimes I_{Z''}) > 0$, then we have $F \cdot L \geq 2(r-1) \cdot L^2$.*

*Proof.* Choose $2(r-1)$ distinct smooth curves $L_1, \ldots, L_{2(r-1)}$ in the complete linear system $|L|$. Choose a set $Z''_i$ of $2(r-1) \cdot L^2$ many distinct points in the open subset

$$L_i - (\bigcup_{j \neq i} L_j)$$

of $L_i$. Let $Z'' = \bigcup_{i=1}^{2(r-1)} Z''_i$. Suppose that $h^0(X, \mathcal{O}_X(F) \otimes I_{Z''}) > 0$. Then, $F$ is effective. If $F$ contains all the curves $L_i$ as its irreducible components, then

$$F \cdot L \geq 2(r-1) \cdot L^2.$$

If $F$ doesn't have $L_i$ as its irreducible component for some $i$, then $F \cap L_i \supset Z''_i$, and

$$F \cdot L = F \cdot L_i \geq \ell(Z''_i) = 2(r-1) \cdot L^2. \blacksquare$$

Now, for $i = 1, \ldots, (r-1)$, we can choose a reduced 0-cycle $Z_i = Z_i' \bigcup Z_i''$ such that $Z_i'$ is chosen as in Lemma 2.5 and $Z_i''$ is chosen as in Lemma 2.6; moreover, we may assume that $Z_1, \ldots, Z_{r-1}$ are disjoint. Put $Z = \bigcup_{i=1}^{r-1} Z_i$, and

$$W = \bigoplus_{i=1}^{(r-1)} [\mathcal{O}_X(L) \otimes I_{Z_i}].$$

Since $h^0(X, \mathcal{O}_X(rL + K_X - c_1) \otimes I_{Z_i' - x}) = 0$ for any $x \in Z_i'$,

$$h^0(X, \mathcal{O}_X(rL + K_X - c_1) \otimes I_{Z_i' \cup Z_i'' - x}) = 0$$

for any $x \in Z_i = Z_i' \bigcup Z_i''$. Hence $Z_i$ satisfies the Cayley-Bacharach property with respect to $|rL + K_X - c_1|$. By Corollary 2.3, there is a bundle $V$ sitting in an extension:

$$0 \to \mathcal{O}_X(c_1 + (1-r)L) \to V \to W \to 0. \tag{2.7}$$

Note that $c_1(V) = c_1$ and that since $Z$ is nonempty, the extension (2.7) is nontrivial.

### 2.3. $L$-Stability of the vector bundle $V$

In the following, we show the $L$-stability of the bundle $V$ constructed above.

**Lemma 2.8.** *The rank-$r$ bundle $V$ in (2.7) is $L$-stable.*

*Proof.* Let $U$ be a proper sub-vector bundle of $V$ such that the quotient $V/U$ is torsion free. Let $U_2$ be the image of $U$ in $W$, and let $U_1$ be the kernel of the surjection $U \to U_2 \to 0$. Then, we have a commutative diagram of morphisms:

$$\begin{array}{ccccccccc}
0 & \to & \mathcal{O}_X(c_1 + (1-r)L) & \to & V & \to & W & \to & 0 \\
& & \uparrow & & \uparrow & & \uparrow & & \\
0 & \to & U_1 & \to & U & \to & U_2 & \to & 0. \\
& & \uparrow & & \uparrow & & \uparrow & & \\
& & 0 & & 0 & & 0 & &
\end{array}$$

Case (a): $U_1 \neq 0$. Then, $c_1(U_1) = (c_1 + (1-r)L) - E_1$ for some effective divisor $E_1$. From $U_2 \hookrightarrow W$, we have $U_2^{**} \hookrightarrow W^{**} = \mathcal{O}_X(L)^{\oplus(r-1)}$; thus,

$$\wedge^{r_2}(U_2^{**}) \hookrightarrow \wedge^{r_2}(\mathcal{O}_X(L)^{\oplus(r-1)}) = \mathcal{O}_X(r_2 L)^{\oplus \binom{r-1}{r_2}}$$

where $r_2$ is the rank of $U_2$. Thus, $c_1(U_2) = r_2 L - E_2$ for some effective divisor $E_2$, and

$$c_1(U) = (c_1 + (1 + r_2 - r)L) - (E_1 + E_2).$$

It follows that $c_1(U) \cdot L \leq (c_1 + (1 + r_2 - r)L) \cdot L$. Therefore,

$$\mu_L(U) = \frac{c_1(U) \cdot L}{(1 + r_2)} \leq \frac{(c_1 + (1 + r_2 - r)L) \cdot L}{(1 + r_2)} < \frac{c_1 \cdot L}{r} = \mu_L(V).$$

Case (b): $U_1 = 0$. Then, $U \hookrightarrow W$; thus, we see that

$$\wedge^{\bar{r}}(U) \hookrightarrow \wedge^{\bar{r}}(W) = \bigoplus_\beta [\mathcal{O}_X(\bar{r}L) \otimes I_{\cup_{i \in \beta} Z_i}]$$

where $\bar{r}$ denotes the rank of $U$ and $\beta$ runs over the set of $\bar{r}$ choices from $(r-1)$ letters. It follows that for some $\beta$ and for some $i \in \beta$, $h^0(X, \mathcal{O}_X(\bar{r}L - c_1(U)) \otimes I_{Z_i}) > 0$. In particular, $h^0(X, \mathcal{O}_X(\bar{r}L - c_1(U)) \otimes I_{Z_i''}) > 0$. In view of Lemma 2.6, we have

$$(\bar{r}L - c_1(U)) \cdot L \geq 2(r-1)L^2 \geq 2\bar{r}L^2.$$

So $c_1(U) \cdot L \leq -\bar{r}L^2 < \bar{r} \cdot (c_1 \cdot L)/r$, and $\mu_L(U) < \mu_L(V)$.

Thus, in both cases, $\mu_L(U) < \mu_L(V)$. Therefore, $V$ is $L$-stable. ∎

In the next lemma, we are going to prove that $h^2(X, \mathrm{ad}(V)) = 0$, that is, the irreducible component of $\mathcal{M}_L(r; c_1, c_2)$ containing $V$ is generically smooth (equivalently, this means that the versal deformation space of $V$ is smooth).

**Lemma 2.9.** *Let $V$ be the rank-$r$ bundle in (2.7). If $rL^2 > K_X \cdot L$, then*
  (i) $\mathrm{Hom}(W, V \otimes \mathcal{O}_X(K_X)) = 0$;
  (ii) $h^2(X, \mathrm{ad}(V)) = 0$.

*Proof.* (i) Let $\beta \in \mathrm{Hom}(W, V \otimes \mathcal{O}_X(K_X))$. Then, $\beta$ induces a map $\beta'$ from $W^{**}$ to $V \otimes \mathcal{O}_X(K_X)$ such that we have commutative diagram of maps:

$$\begin{array}{ccc} W & \hookrightarrow & W^{**} = \mathcal{O}_X(L)^{\oplus(r-1)} \\ {\scriptstyle \beta}\downarrow & \swarrow {\scriptstyle \beta'} & \\ V \otimes \mathcal{O}_X(K_X). & & \end{array}$$

To show that $\beta = 0$, it suffices to show that $H^0(X, V \otimes \mathcal{O}_X(K_X - L)) = 0$.

Since $c_1 \cdot L \leq 0$ and $K_X \cdot L < rL^2$, $(c_1 - rL + K_X) \cdot L < 0$. Thus,

$$H^0(X, \mathcal{O}_X(c_1 - rL + K_X)) = 0.$$

By our choice of the 0-cycles $Z_i'$, $H^0(X, \mathcal{O}_X(K_X) \otimes I_{Z_i'}) = 0$. Thus,

$$H^0(X, W \otimes \mathcal{O}_X(K_X - L)) = 0.$$

Now, tensoring (2.7) by $\mathcal{O}_X(K_X - L)$ and taking cohomology, we see that

$$H^0(X, V \otimes \mathcal{O}_X(K_X - L)) = 0.$$

(ii) We follow the argument as in the proof of Lemma 4.5.4 in [19]. By the Serre duality, we have $H^2(X, \mathrm{ad}(V)) \cong H^0(X, \mathrm{ad}(V) \otimes \mathcal{O}_X(K_X))$. Let

$$\phi \in H^0(X, \mathrm{ad}(V) \otimes \mathcal{O}_X(K_X)) \subseteq H^0(X, \mathcal{E}nd(V) \otimes \mathcal{O}_X(K_X)).$$

Then, we obtain a map $\phi$ from $V$ to $V \otimes \mathcal{O}_X(K_X)$. Consider the diagram:

$$0 \to \mathcal{O}_X(c_1 + (1-r)L) \xrightarrow{\theta} V \xrightarrow{\rho} W \to 0 \qquad (2.10)$$

$$\downarrow \phi$$

$$0 \to \mathcal{O}_X(c_1 + (1-r)L + K_X) \xrightarrow{\theta'} V \otimes \mathcal{O}_X(K_X) \xrightarrow{\rho'} W \otimes \mathcal{O}_X(K_X) \to 0. \qquad (2.11)$$

By our choice of the 0-cycles $Z'_i$, $H^0(X, \mathcal{O}_X(rL - c_1 + K_X) \otimes I_{Z'_i}) = 0$. Thus,

$$\mathrm{Hom}(\mathcal{O}_X(c_1 + (1-r)L), W \otimes \mathcal{O}_X(K_X)) = 0,$$

so $\rho' \circ \phi \circ \theta = 0$. Applying $\mathrm{Hom}(\mathcal{O}_X(c_1 + (1-r)L), \cdot)$ to (2.11), we obtain:

$$0 \to H^0(X, \mathcal{O}_X(K_X)) \xrightarrow{\lambda} \mathrm{Hom}(\mathcal{O}_X(c_1 + (1-r)L), V \otimes \mathcal{O}_X(K_X))$$
$$\xrightarrow{\rho' \circ} \mathrm{Hom}(\mathcal{O}_X(c_1 + (1-r)L), W \otimes \mathcal{O}_X(K_X)) = 0.$$

It follows that there exists $\tau \in H^0(X, \mathcal{O}_X(K_X))$ such that

$$\phi \circ \theta = \lambda(\tau) = (\tau \otimes \mathrm{Id}_V) \circ \theta$$

where $\mathrm{Id}_V$ is the identity morphism in $\mathrm{End}(V)$. Thus, $(\phi - \tau \otimes \mathrm{Id}_V) \circ \theta = 0$. Applying $\mathrm{Hom}(\cdot, V \otimes \mathcal{O}_X(K_X))$ to (2.10), we get an exact sequence:

$$\mathrm{Hom}(W, V \otimes \mathcal{O}_X(K_X)) \to H^0(X, \mathcal{E}nd(V) \otimes \mathcal{O}_X(K_X))$$
$$\xrightarrow{\circ \theta} \mathrm{Hom}(\mathcal{O}_X(c_1 + (1-r)L), V \otimes \mathcal{O}_X(K_X)).$$

From (i), we conclude that $(\phi - \tau \otimes \mathrm{Id}_V) = 0$. Since $0 = \mathrm{Tr}(\phi) = \tau$, $\phi = 0$. Hence,

$$h^2(X, \mathrm{ad}(V)) = 0. \blacksquare$$

Finally, we state and prove the main result in this section.

**Theorem 2.12.** *For any ample divisor $L$ on $X$, there exists a constant $\alpha$ depending only on $X$, $r$, $c_1$ and $L$ such that for any $c_2 \geq \alpha$, there exists an $L$-stable rank-$r$ bundle $V$ with Chern classes $c_1$ and $c_2$. Moreover, $h^2(X, \mathrm{ad}(V)) = 0$.*

*Proof.* We may re-scale the ample divisor $L$ such that $L$ is very ample and that $rL^2 > K_X \cdot L$. Note that $c_1(W) = (r-1)L$ and $c_2(W) = \ell(Z) + (r-1)(r-2)/2 \cdot L^2$. From the exact sequence (2.7), we see that $c_1(V) = c_1$ and

$$c_2(V) = \ell(Z) + (r-1)(c_1 \cdot L) - r(r-1)/2 \cdot L^2.$$

By the construction of the 0-cycle $Z$, we get

$$\ell(Z) = \sum_{i=1}^{(r-1)} [\ell(Z'_i) + \ell(Z''_i)]$$
$$\geq (r-1)[1 + \max(p_g, h^0(X, \mathcal{O}_X(rL - c_1 + K_X))) + 4(r-1)^2 \cdot L^2].$$

Let $\alpha$ be the integer:

$$(r-1)[1 + \max(p_g, h^0(X, \mathcal{O}_X(rL - c_1 + K_X))) + 4(r-1)^2 \cdot L^2]$$
$$+ (r-1)(c_1 \cdot L) - r(r-1)/2 \cdot L^2.$$

Then, $\alpha$ depends only on $X$, $r$, $c_1$ and $L$. By Lemma 2.8, for any $c_2 \geq \alpha$, there exists an $L$-stable rank-$r$ bundle $V$ with Chern classes $c_1$ and $c_2$.

Moreover, since $rL^2 > K_X \cdot L$, $h^2(X, \mathrm{ad}(V)) = 0$ by Lemma 2.9 (ii). ∎

*Remark 2.13.* In [2], Artamkin showed that $\mathcal{M}_L(r; 0, c_2)$ is nonempty whenever

$$c_2 > (r+1) \cdot \max(1, p_g);$$

in particular, when we only consider the case of $c_1 = 0$, the lower bound of the integer $c_2$ does not depend on the ample divisor $L$. By contrast, the constant $\alpha$ in Theorem 2.12 depends on $L$. In fact, if we want a universal lower bound of $c_2$ for all $c_1$, this bound must depend on the ample divisor $L$. We shall see this fact from Theorem 3.1 in the next section that on a ruled surface, there exists a divisor $c_1$ such that for any integer $c_2$, we can find an ample divisor $L$ with $\mathcal{M}_L(r; c_1, c_2)$ being empty.

**3. Restriction of a stable bundle on a ruled surface to the generic fiber**

From now on, we study stable bundles on a ruled surface $X$. Our first goal in this section is to show that if $0 < (c_1 \cdot f) < r$ and if $r_L \gg 0$, then $\mathcal{M}_L(r; c_1, c_2)$ is empty.

**Theorem 3.1.** *Let $0 < (c_1 \cdot f) < r$. Then, there exists a constant $r_0$ depending only on $X, r, c_1$ and $c_2$ such that $\mathcal{M}_L(r; c_1, c_2)$ is empty whenever $r_L > r_0$.*

*Proof.* Assume that $V \in \mathcal{M}_L(r; c_1, c_2)$. Let $c_1 = (a\sigma + \mathbf{b}f)$; then, $0 < a < r$. For any divisor $\mathbf{k}$ on $C$, we see that $c_1(V \otimes \mathcal{O}_X(-\sigma + \mathbf{k}f)) = (a-r)\sigma + (\mathbf{b} + r\mathbf{k})f$ and that

$$c_2(V \otimes \mathcal{O}_X(-\sigma + \mathbf{k}f)) = c_2 + (r-1)(a\sigma + \mathbf{b}f) \cdot (-\sigma + \mathbf{k}f) + \frac{r(r-1)}{2} \cdot (-\sigma + \mathbf{k}f)^2.$$

By the Riemann-Roch formula, we conclude the following:

$$\chi(V \otimes \mathcal{O}_X(-\sigma + \mathbf{k}f)) = a \cdot k + a \cdot (b + 1 - g_C) - c_2 - \frac{e(a^2 - a)}{2}.$$

Let $k = g_C - b + [c_2/a + e(a-1)/2] + 1$. Then, $\chi(V \otimes \mathcal{O}_X(-\sigma + \mathbf{k}f)) > 0$. Thus, $h^i(X, V \otimes \mathcal{O}_X(-\sigma + \mathbf{k}f)) > 0$ where $i = 0$ or $2$. On the other hand, put

$$r_0 = \max\{e + \frac{kr + b}{r - a}, \quad e - \frac{2r\chi(\mathcal{O}_X) + er + kr + b}{r + a}\}.$$

Then, $r_0$ is a number depending only on $X, r, c_1$ and $c_2$. If $h^0(X, V \otimes \mathcal{O}_X(-\sigma + \mathbf{k}f)) > 0$, then there exists an injective map $\mathcal{O}_C(\sigma - \mathbf{k}f) \hookrightarrow V$. By stability of $V$, we see that $(\sigma - \mathbf{k}f) \cdot L < (a\sigma + bf) \cdot L/r$. By direct calculations, we get

$$r_L < e + \frac{kr + b}{r - a};$$

but this contradicts with the choice of the numbers $r_0$ and $r_L$.

If $h^2(X, V \otimes \mathcal{O}_X(-\sigma + \mathbf{k}f)) > 0$, then $h^0(X, V^* \otimes \mathcal{O}_X(K_X + \sigma - \mathbf{k}f)) > 0$. Hence, there is a nonzero map $V \to \mathcal{O}_X(K_X - \sigma + \mathbf{k}f)$ which can be extended to

$$V \to \mathcal{O}_X(K_X + \sigma - \mathbf{k}f) \otimes \mathcal{O}_X(-E) \otimes I_Z \to 0$$

for some effective divisor $E$. By the stability of $V$, we must have

$$c_1(V) \cdot L/r < K_X \cdot L + (\sigma - \mathbf{k}f) \cdot L - E \cdot L \le K_X \cdot L + (\sigma - \mathbf{k}f) \cdot L.$$

By a straightforward calculation, we obtain that

$$r_L \le e - \frac{2r\chi(\mathcal{O}_X) + er + kr + b}{r + a};$$

again, this contradicts with our choices of $r_0$ and $r_L$.

Therefore, if $r_L > r_0$, the moduli space $\mathcal{M}_L(r; c_1, c_2)$ is empty. ∎

*Remark 3.2.* Theorem 3.1 only says that for a fixed $c_1$ with $0 < c_1 \cdot f < r$ and for a fixed $c_2$, the moduli space $\mathcal{M}_L(r; c_1, c_2)$ is empty for some special ample divisor $L$ (e.g., when $r_L > r_0$). For other ample divisor $L'$, $\mathcal{M}_{L'}(r; c_1, c_2)$ can be nonempty (see [21] when $r = 2$); we will discuss this issue in other places.

In view of Theorem 3.1, our next goal is to study the moduli space $\mathcal{M}_L(r; \mathbf{t}f, c_2)$ where $-r < t \leq 0$. Let $V \in \mathcal{M}_L(r; \mathbf{t}f, c_2)$ where $L$ is of the form $(\sigma + \mathbf{r_L}f)$ with

$$r_L \geq \max\{e/2 - \chi(\mathcal{O}_X) + r(g_C + |c_2|) + 1,\ 2|e| + r(g_C + |c_2|)\}. \tag{3.3}$$

We want to show that the restriction of the stable bundle $V$ to the generic fiber is trivial. To start with, we prove the following technical lemma.

**Lemma 3.4.** *Let $U$ be a rank-$s$ bundle with an injection $U \hookrightarrow V$.*
  (i) *For any divisor $\mathbf{d}$ with $d \geq -r(g_C + |c_2|) - 1$, $h^2(X, U^* \otimes \mathcal{O}_X(\mathbf{d}f)) = 0$;*
  (ii) *If $c_1(U) = -\mathbf{a}f$ with $0 < a \leq (r-s)(g_C + |c_2|)$ and $c_2(U) \leq c_2$, then $U$ sits in*

$$0 \to U_1 \to U \to \mathcal{O}_X(\mathbf{n}f) \otimes I_Z \to 0$$

*where $U_1$ is a rank-$(s-1)$ bundle with an injection $U_1 \hookrightarrow V$; moreover, $c_1(U_1) = -(\mathbf{a} + \mathbf{n})f$ with $0 < (a + n) \leq (r - s + 1)(g_C + |c_2|)$, and $c_2(U_1) \leq c_2$.*

*Proof.* (i) By the Serre duality, $h^2(X, U^* \otimes \mathcal{O}_X(\mathbf{d}f)) = h^0(X, U \otimes \mathcal{O}_X(K_X - \mathbf{d}f))$. If

$$h^0(X, U \otimes \mathcal{O}_X(K_X - \mathbf{d}f)) > 0,$$

then we have $\mathcal{O}_X(\mathbf{d}f - K_X) \hookrightarrow U \hookrightarrow V$; by the stability of $V$, we obtain that

$$(\mathbf{d}f - K_X) \cdot L < \frac{tf \cdot L}{r} \leq 0.$$

On the other hand, we have $(\mathbf{d}f - K_X) \cdot L = d - 2(e/2 - \chi(\mathcal{O}_X)) + 2r_L \geq 0$ in view of the assumption (3.3); but this is a contradiction.

(ii) By the Riemann-Roch formula, one checks that

$$\chi(U^* \otimes \mathcal{O}_X(\mathbf{k}f)) = s \cdot k + s \cdot \chi(\mathcal{O}_X) + a - c_2(U) \geq s \cdot k + s \cdot \chi(\mathcal{O}_X) + a - c_2.$$

Let $k = g_C + [(c_2 - a)/s]$. Then, $\chi(U^* \otimes \mathcal{O}_X(\mathbf{k}f)) > 0$. Since

$$k \geq g_C + \frac{c_2 - (r-s)(g_C + |c_2|)}{s} - 1 \geq -r(g_C + |c_2|) - 1,$$

$h^0(X, U^* \otimes \mathcal{O}_X(\mathbf{k}f)) > 0$ by (i); thus, there is an exact sequence:

$$0 \to U_1 \to U \to \mathcal{O}_X(\mathbf{k}f - E) \otimes I_Z \to 0$$

where $E$ is effective and $Z$ is a 0-cycle. Since $U/U_1$ is torsion-free, $U_1$ is a bundle. Let $E \equiv (\lambda \sigma + \mu f)$. Then, $\lambda \geq 0$; moreover, $\mu \geq 0$ when $e \geq 0$, and $\mu \geq \lambda e/2$ when $e < 0$. We claim that $\lambda = 0$: otherwise, $\lambda \geq 1$; then

$$c_1(U_1) \cdot L = (\lambda \sigma + (\mu - a - k)f) \cdot L$$
$$= \lambda(r_L - e) + \mu - a - k$$
$$\geq (r_L - e) - |e| - a - k.$$

But

$$a + k \leq (r - s)(g_C + |c_2|) + g_C + [(c_2 - a)/s]$$
$$\leq (r - s)(g_C + |c_2|) + g_C + |c_2|$$
$$= (r - s + 1)(g_C + |c_2|)$$
$$\leq r(g_C + |c_2|).$$

So $c_1(U_1) \cdot L \geq r_L - 2|e| - r(g_C + |c_2|) \geq 0$ by our assumption about $r_L$; but this contradicts with the stability of $V$. Therefore, $E$ is supported in the fibers of the ruling, and $U$ sits in the desired exact sequence; moreover, $c_2(U_1) \leq c_2(U) \leq c_2$. Note that $c_1(U_1) = -(\mathbf{a} + \mathbf{n})f$ and that $(a + n) \leq (a + k) \leq (r - s + 1)(g_C + |c_2|)$. By the stability of $V$, $-(a + n)/(s - 1) < -t/r \leq 0$. Thus, $(a + n) > 0$. ∎

**Theorem 3.5.** *Let $V \in \mathcal{M}_L(r; \mathbf{t}f, c_2)$ where $-r < t \leq 0$ and $L$ satisfies (3.3). Then,*

$$V|_{\mathbf{P}_K^1} = \mathcal{O}_{\mathbf{P}_K^1}^{\oplus r}.$$

*Proof.* By Lemma 3.4 (ii) and by induction on the rank of subbundles of $V$, we conclude that there exists a flag of subbundles of $V$: $V_1 \subset V_2 \subset \ldots \subset V_{r-1} \subset V_r = V$ such that $\mathrm{rank}(V_i) = i$, $c_2(V_i) \leq c_2$, $c_1(V_i) = -\mathbf{b}_i f$ with $0 < b_i \leq r(g_C + |c_2|)$ for $i < r$, and $V_i/V_{i-1} = \mathcal{O}_X((\mathbf{b}_{i-1} - \mathbf{b}_i)f) \otimes I_{Z_i}$ where $Z_i$ is a 0-cycle. Hence $V|_{\mathbf{P}_K^1} = \mathcal{O}_{\mathbf{P}_K^1}^{\oplus r}$. ∎

Next, we prove the following simple observation.

**Lemma 3.6.** *If the moduli space $\mathcal{M}_L(r; \mathbf{t}f, c_2)$ is nonempty, then it is smooth with dimension $2rc_2 - (r^2 - 1)(1 - g_C)$; in particular, $c_2 \geq (1 - g_C)(r^2 - 1)/(2r)$.*

*Proof.* Since $L$ satisfies (3.3), $K_X \cdot L \leq 0$. By a well-known result of Maruyama, $\mathcal{M}_L(r; \mathbf{t}f, c_2)$ is smooth with the expected dimension $2rc_2 - (r^2 - 1)(1 - g_C)$. ∎

We notice that the ample divisor $L$ in Theorem 3.5 depends on the integer $c_2$ (that is, the condition (3.3)). However, in our existence result Theorem 2.12, the integer $c_2$ has to be bigger than some constant depending on $L$. Thus, Theorem 2.12 can not apply to the present situation to guarantee the nonemptiness of the moduli space $\mathcal{M}_L(r; \mathbf{t}f, c_2)$. The following result deals with this problem.

**Proposition 3.7.** *Let $r \geq 2$, $-r < t \leq 0$, and $L = (\sigma + \mathbf{r_L}f)$ with $r_L \geq (|e| + 2r - 2)$. If $c_2 \geq 2(r - 1)$, then the moduli space $\mathcal{M}_L(r; \mathbf{t}f, c_2)$ is nonempty.*

We omit the proof since it is a slight modification of the proof of Theorem 2.12 (replacing the $L$ in $W$ by $f$). It seems to us that a stronger result should hold, that is, if $c_2 \geq (r + t)$, then $\mathcal{M}_L(r; \mathbf{t}f, c_2)$ is nonempty (see Theorem 5.4 (iii)).

**4. Generic bundles in $\mathcal{M}_L(r; tf, c_2)$ on a rational ruled surface**

From now on, $X$ will be a rational ruled surface. In this section, we will study the structure of a generic bundle in $\mathcal{M}_L(r; tf, c_2)$ where $L$ satisfies (3.3) and $-r < t \leq 0$.

**4.1.** Exact sequences associated to a bundle $V$ in $\mathcal{M}_L(r; tf, c_2)$

In this subsection, we will construct $(r - 1)$ exact sequences for each vector bundle in the moduli space $\mathcal{M}_L(r; tf, c_2)$. We begin with two lemmas.

**Lemma 4.1.** *Let $U$ be a rank-$i$ bundle with $c_1(U) = af$ and $U|_{\mathbf{P}^1_K} = \mathcal{O}_{\mathbf{P}^1_K}^{\oplus i}$. Then,*
  (i) *$\pi_* U$ is a rank-$i$ bundle on $\mathbf{P}^1$;*
  (ii) *$\deg c_1(\pi_* U) \geq (a - c_2(U))$.*

*Proof.* (i) Note that $\pi_* U$ is always torsion-free. Thus, $\pi_* U$ is a vector bundle. Since $U|_{\mathbf{P}^1_K}$ is equal to $\mathcal{O}_{\mathbf{P}^1_K}^{\oplus i}$, the rank of $\pi_* U$ is equal to $i$.

(ii) Since $U|_{\mathbf{P}^1_K} = \mathcal{O}_{\mathbf{P}^1_K}^{\oplus i}$, $R^1\pi_* U$ is a torsion sheaf supported in some points; thus, $\deg c_1(R^1\pi_* U) \geq 0$. By the Grothendieck-Riemann-Roch formula (see p.436 in [12]),

$$\text{ch}(\pi_* U) - \text{ch}(R^1\pi_* U) = \pi_*(\text{ch}(U) \cdot \text{td}(T_\pi)) = i + (a - c_2(U)) \cdot [pt]$$

where $T_\pi$ is the relative tangent bundle, $\mathrm{td}(T_\pi) = 1 + (\sigma - e/2 \cdot f)$, and $[pt]$ stands for the class determined by a point. Therefore,

$$\deg c_1(\pi_* U) = \deg c_1(R^1\pi_* U) + (a - c_2(U)) \geq (a - c_2(U)).\blacksquare$$

**Lemma 4.2.** *Let $U$ be a rank-$i$ bundle with $c_1(U) = af$ and $U|_{\mathbf{P}^1_K} = \mathcal{O}^{\oplus i}_{\mathbf{P}^1_K}$. If*

$$\pi_* U = \mathcal{O}_{\mathbf{P}^1}(-n)^{\oplus j} \oplus \mathcal{O}_{\mathbf{P}^1}(-n_1) \oplus \ldots \oplus \mathcal{O}_{\mathbf{P}^1}(-n_{i-j}) \tag{4.3}$$

*where $1 \leq j \leq i$ and $n < n_1 \leq \ldots \leq n_{i-j}$, then*
(i) $in + (i - j) \leq (c_2(U) - a)$;
(ii) $in - h^0(\mathbf{P}^1, \pi_* U \otimes \mathcal{O}_{\mathbf{P}^1}(n)) \leq (c_2(U) - a) - i$;
(iii) *the bundle $U$ sits in an exact sequence of the form:*

$$0 \to \mathcal{O}_X(-nf) \to U \to W \to 0 \tag{4.4}$$

*where $W$ is a torsion-free rank-$(i-1)$ sheaf with $W|_{\mathbf{P}^1_K} = (W^{**})|_{\mathbf{P}^1_K} = \mathcal{O}^{\oplus(i-1)}_{\mathbf{P}^1_K}$.*

*Proof.* (i) Since $n < n_1 \leq \ldots \leq n_{i-j}$, by Lemma 4.1 (ii), we have

$$(c_2(U) - a) \geq -\deg c_1(\pi_* U) = jn + \sum_{k=1}^{i-j} n_k \geq in + (i - j).$$

(ii) Note that $h^0(\mathbf{P}^1, \pi_* U \otimes \mathcal{O}_{\mathbf{P}^1}(n)) = j$. Therefore, by (i),

$$in - h^0(\mathbf{P}^1, \pi_* U \otimes \mathcal{O}_{\mathbf{P}^1}(n)) \leq [(c_2(U) - a) - (i - j)] - j = (c_2(U) - a) - i.$$

(iii) Since there is a natural injection $\pi^*(\pi_* U) \hookrightarrow U$, we have

$$\mathcal{O}_X(-nf) \hookrightarrow U.$$

We claim that the quotient $W = U/\mathcal{O}_X(-nf)$ is torsion free: otherwise, we have

$$\mathcal{O}_X(-nf) \hookrightarrow \mathcal{O}_X(-nf + D) \hookrightarrow U \tag{4.5}$$

where $D$ is some nontrivial effective divisor; since $U|_{\mathbf{P}^1_K}$ is equal to $\mathcal{O}^{\oplus i}_{\mathbf{P}^1_K}$, $D$ is supported in the fibers of $\pi$; put $D = df$ where $d > 0$; applying $\pi_*$ to (4.5), we obtain

$$\mathcal{O}_{\mathbf{P}^1}(-n) \hookrightarrow \mathcal{O}_{\mathbf{P}^1}(-n + d) \hookrightarrow \pi_* U;$$

but this is impossible in view of the assumption (4.3).

Thus, we have the exact sequence (4.4). Since $W$ is torsion-free, $W$ is locally free outside possibly finitely many points. Restricting (4.4) to $\mathbf{P}_K^1$, we see that

$$0 \to \mathcal{O}_{\mathbf{P}_K^1} \to \mathcal{O}_{\mathbf{P}_K^1}^{\oplus i} \to W|_{\mathbf{P}_K^1} = W^{**}|_{\mathbf{P}_K^1} \to 0.$$

Since $c_1(W) = (a+n)f$, we conclude that $W|_{\mathbf{P}_K^1} = W^{**}|_{\mathbf{P}_K^1} = \mathcal{O}_{\mathbf{P}_K^1}^{\oplus(i-1)}$. ∎

**Proposition 4.6.** *Let $V \in \mathcal{M}_L(r; tf, c_2)$ where $L$ satisfies the condition (3.3) and $-r < t \leq 0$. Then, there exist $(r-1)$ exact sequences:*

$$0 \to \mathcal{O}_X(-n_i f) \to V_i^{**} \to V_{i-1} \to 0 \tag{4.7}$$

*where $i = r, \ldots, 2$, $V_r = V$, and $V_i$ is a torsion-free rank-$i$ sheaf such that*

(i) $\pi_*(V_i^{**}) = \mathcal{O}_{\mathbf{P}^1}(-n_i)^{\oplus j_i} \oplus \mathcal{O}_{\mathbf{P}^1}(-n_{i,1}) \oplus \ldots \oplus \mathcal{O}_{\mathbf{P}^1}(-n_{i,i-j_i})$ *with $n_i < n_{i,k}$;*

(ii) $(V_i)|_{\mathbf{P}_K^1} = (V_i^{**})|_{\mathbf{P}_K^1} = \mathcal{O}_{\mathbf{P}_K^1}^{\oplus i}$;

(iii) $in_i + (i - j_i) \leq (c_2(V_i^{**}) - t - \sum_{k=i+1}^{r} n_k)$;

(iv) $in_i - h^0(\mathbf{P}^1, \pi_*(V_i^{**}) \otimes \mathcal{O}_{\mathbf{P}^1}(n_i)) \leq (c_2(V_i^{**}) - t - \sum_{k=i+1}^{r} n_k) - i$.

*Proof.* By Theorem 3.5, $V|_{\mathbf{P}_K^1} = \mathcal{O}_{\mathbf{P}_K^1}^{\oplus r}$. Now, the exact sequences (4.7) and the properties (i) and (ii) follow from induction and Lemma 4.2 (iii). Note that

$$c_1(V_i^{**}) = c_1(V_i) = c_1(V_{i+1}^{**}) + n_{i+1}f = \left(t + \sum_{k=i+1}^{r} n_k\right) f.$$

Therefore, the properties (iii) and (iv) follow from Lemma 4.2 (i) and (ii). ∎

**4.2. The number of moduli of $V_i$ and $V_i^{**}$**

In this subsection, we estimate the number of moduli of $V_i$ and $V_i^{**}$. These estimations will be used in the next subsection to study generic bundles in the moduli space $\mathcal{M}_L(r; tf, c_2)$ where $L$ satisfies the condition (3.3) and $-r < t \leq 0$. To begin with, we collect some properties satisfied by the sheaf $V_i$.

**Lemma 4.8.** (i) *For each $i$, there exists a canonical exact sequence*

$$0 \to V_i \to V_i^{**} \to Q_i \to 0 \tag{4.9}$$

*where $Q_i$ is a torsion sheaf supported on finitely many points in $X$;*

*(ii) $\dim \operatorname{Hom}(V_i, \mathcal{O}_X(-n_{i+1}f)) + 1 \leq \dim \operatorname{Aut}(V_{i+1}^{**})$;*

*(iii) $\operatorname{Ext}^2(V_i, \mathcal{O}_X(-n_{i+1}f)) = 0$;*

*(iv) $-\chi(V_i, \mathcal{O}_X(-n_{i+1}f)) = c_2(V_i) + (t + \sum_{k=i+1}^{r} n_k) + i \cdot n_{i+1} - i.$*

*Proof.* (i) This is a standard fact. The torsion sheaf $Q_i$ is supported on those points where $V_i$ is not locally free.

(ii) Applying the functor $\operatorname{Hom}(V_{i+1}^{**}, \cdot)$ to the exact sequence (4.10), we have

$$0 \to \operatorname{Hom}(V_{i+1}^{**}, \mathcal{O}_X(-n_{i+1}f)) \to \operatorname{End}(V_{i+1}^{**}) \xrightarrow{\psi_{i+1}} \operatorname{Hom}(V_{i+1}^{**}, V_i)$$

where $\psi_{i+1}(\operatorname{Id}) = p_i$ for the identity endormorphism Id in $\operatorname{End}(V_{i+1}^{**})$. Thus,

$$\dim \operatorname{Aut}(V_{i+1}^{**}) = \dim \operatorname{End}(V_{i+1}^{**}) \geq 1 + \dim \operatorname{Hom}(V_{i+1}^{**}, \mathcal{O}_X(-n_{i+1}f)).$$

Similarly, applying the functor $\operatorname{Hom}(\cdot, \mathcal{O}_X(-n_{i+1}f))$ to (4.10), we obtain

$$0 \to \operatorname{Hom}(V_i, \mathcal{O}_X(-n_{i+1}f)) \to \operatorname{Hom}(V_{i+1}^{**}, \mathcal{O}_X(-n_{i+1}f));$$

thus, $\dim \operatorname{Hom}(V_{i+1}^{**}, \mathcal{O}_X(-n_{i+1}f)) \geq \dim \operatorname{Hom}(V_i, \mathcal{O}_X(-n_{i+1}f))$. Hence,

$$\dim \operatorname{Aut}(V_{i+1}^{**}) \geq 1 + \dim \operatorname{Hom}(V_i, \mathcal{O}_X(-n_{i+1}f)).$$

(iii) Since $\mathcal{O}_X(K_X + n_{i+1}f)|_{\mathbf{P}_K^1} = \mathcal{O}_{\mathbf{P}_K^1}(-2)$ and $(V_i)|_{\mathbf{P}_K^1} = \mathcal{O}_{\mathbf{P}_K^1}^{\oplus i}$, we see that $H^0(X, V_i \otimes \mathcal{O}_X(K_X + n_{i+1}f)) = 0$. By the Serre duality,

$$\operatorname{Ext}^2(V_i, \mathcal{O}_X(-n_{i+1}f)) \cong H^0(X, V_i \otimes \mathcal{O}_X(K_X + n_{i+1}f)) = 0.$$

(iv) Recall that by definition, $\chi(\mathcal{F}_1, \mathcal{F}_2) = \sum_{i=0}^{2} (-1)^i \dim \operatorname{Ext}^i(\mathcal{F}_1, \mathcal{F}_2)$ for two sheaves $\mathcal{F}_1$ and $\mathcal{F}_2$ on $X$. Let $\operatorname{td}(X)$ be the Todd class of $X$, and let $\operatorname{ch}(\mathcal{F})$ be the Chern character of a sheaf $\mathcal{F}$. Then, we have the formula:

$$\chi(\mathcal{F}_1, \mathcal{F}_2) = (\operatorname{ch}(\mathcal{F}_1)^* \cdot \operatorname{ch}(\mathcal{F}_2) \cdot \operatorname{td}(X))_4$$

where $*$ acts on $H^{2i}(X, \mathbf{Z})$ by multiplication of $(-1)^i$. Thus, we obtain

$$-\chi(V_i, \mathcal{O}_X(-n_{i+1}f)) = c_2(V_i) - \frac{1}{2}(K_X \cdot c_1(V_i)) + i \cdot n_{i+1} - i.$$

Since $c_1(V_i) = (t + \sum_{k=i+1}^{r} n_k)f$, the conclusion follows immediately. ∎

Next, for convenience, we introduce some notations.

*Notation 4.10.* (i) Let $\ell_i = h^0(X, Q_i)$ for $i = 1, \ldots, r-1$;

(ii) Let $\delta_i = [\#(\text{moduli of } V_i) - \dim \text{Aut}(V_i)]$ for $i = 1, \ldots, r-1$;

(iii) Let $\delta_i^{**} = [\#(\text{moduli of } V_i^{**}) - \dim \text{Aut}(V_i^{**})]$ for $i = 1, \ldots, r$.

Now, we estimate the number of moduli of $Q_i$, $V_i$ and $V_i^{**}$.

**Lemma 4.11.** (i) $\#(\text{moduli of } Q_i) - \dim \text{Aut}(Q_i) \leq \ell_i$;

(ii) $\delta_i \leq \delta_i^{**} + (i+1)\ell_i$;

(iii) $\delta_i^{**} \leq \delta_{i-1} - \chi(V_{i-1}, \mathcal{O}_X(-n_i f)) - h^0(X, V_i^{**} \otimes \mathcal{O}_X(n_i f))$.

*Proof.* (i) From (4.7), we have an exact sequence

$$0 \to \mathcal{O}_X(-n_{i+1}f) \to V_{i+1}^{**} \to V_i \to 0. \tag{4.12}$$

Applying $\text{Hom}(\cdot, Q_i)$ to (4.12), we obtain

$$\dim \text{Hom}(V_i, Q_i) \leq \dim \text{Hom}(V_{i+1}^{**}, Q_i) = (i+1)\ell_i.$$

Applying $\text{Hom}(\cdot, Q_i)$ to (4.9), we get

$$0 \to \text{Hom}(Q_i, Q_i) \to \text{Hom}(V_i^{**}, Q_i) \to \text{Hom}(V_i, Q_i)$$
$$\to \text{Ext}^1(Q_i, Q_i) \to \text{Ext}^1(V_i^{**}, Q_i) = 0.$$

It follows that

$$\dim \text{Ext}^1(Q_i, Q_i) - \dim \text{Hom}(Q_i, Q_i) = \dim \text{Hom}(V_i, Q_i) - \dim \text{Hom}(V_i^{**}, Q_i)$$
$$\leq (i+1)\ell_i - i\ell_i = \ell_i.$$

Since $\#(\text{moduli of } Q_i) \leq \dim \text{Ext}^1(Q_i, Q_i)$ and $\dim \text{Aut}(Q_i) = \dim \text{Hom}(Q_i, Q_i)$,

$$\#(\text{moduli of } Q_i) - \dim \text{Aut}(Q_i) \leq \ell_i. \tag{4.13}$$

(ii) From the exact sequence (4.9), we see that

$$\#(\text{moduli of } V_i) \leq \#(\text{moduli of } V_i^{**}) + \#(\text{moduli of } Q_i) + \dim \text{Hom}(V_i^{**}, Q_i)$$
$$- \dim \text{Aut}(V_i^{**}) - \dim \text{Aut}(Q_i) + 1$$
$$= \delta_i^{**} + [\#(\text{moduli of } Q_i) - \dim \text{Aut}(Q_i)]$$
$$+ \dim \text{Hom}(V_i^{**}, Q_i) + 1$$
$$\leq \delta_i^{**} + \ell_i + i\ell_i + 1$$
$$= \delta_i^{**} + (i+1)\ell_i + 1.$$

Since dim $\text{Aut}(V_i) \geq 1$, we obtain that $\delta_i \leq \delta_i^{**} + (i+1)l_i$.

(iii) Similarly, from the exact sequence (4.7), we have

$$\#(\text{moduli of } V_i^{**}) \leq \#(\text{moduli of } V_{i-1}) + \dim \text{Ext}^1(V_{i-1}, \mathcal{O}_X(-n_i f))$$
$$- \dim \text{Hom}(\mathcal{O}_X(-n_i f), V_i^{**}) - \dim \text{Aut}(V_{i-1}) + 1$$
$$= \delta_{i-1} + \dim \text{Ext}^1(V_{i-1}, \mathcal{O}_X(-n_i f)) - h^0(X, V_i^{**} \otimes \mathcal{O}_X(n_i f)) + 1$$
$$= \delta_{i-1} - \chi(V_{i-1}, \mathcal{O}_X(-n_i f)) - h^0(X, V_i^{**} \otimes \mathcal{O}_X(n_i f))$$
$$+ 1 + \dim \text{Hom}(V_{i-1}, \mathcal{O}_X(-n_i f))$$

where we have used Lemma 4.8 (iii) in the last equality. By Lemma 4.8 (ii),

$$\delta_i^{**} \leq \delta_{i-1} - \chi(V_{i-1}, \mathcal{O}_X(-n_i f)) - h^0(X, V_i^{**} \otimes \mathcal{O}_X(n_i f)). \blacksquare$$

**Proposition 4.14.** $\delta_i^{**} \leq \delta_{i-1}^{**} + 2(c_2 - \sum_{k=i}^{r-1} l_k) - (2i-1) + il_{i-1}$.

*Proof.* By Lemma 4.8 (iv) and Proposition 4.6 (iv), we have

$$-\chi(V_{i-1}, \mathcal{O}_X(-n_i f)) = c_2(V_{i-1}) + (t + \sum_{k=i}^{r} n_k) + (i-1)n_i - (i-1)$$
$$= c_2(V_i^{**}) + (t + \sum_{k=i+1}^{r} n_k + in_i) - (i-1)$$
$$\leq c_2(V_i^{**}) + [c_2(V_i^{**}) + h^0(\mathbf{P}^1, \pi_*(V_i^{**}) \otimes \mathcal{O}_{\mathbf{P}^1}(n_i))] - (2i-1)$$
$$= 2c_2(V_i^{**}) + h^0(X, V_i^{**} \otimes \mathcal{O}_X(n_i f)) - (2i-1)$$
$$= 2(c_2 - \sum_{k=i}^{r-1} l_k) + h^0(X, V_i^{**} \otimes \mathcal{O}_X(n_i f)) - (2i-1).$$

Therefore, by Lemma 4.11 (ii) and (iii), we conclude that

$$\delta_i^{**} \leq \delta_{i-1}^{**} - \chi(V_{i-1}, \mathcal{O}_X(-n_i f)) - h^0(X, V_i^{**} \otimes \mathcal{O}_X(n_i f)) + il_{i-1}$$
$$\leq \delta_{i-1}^{**} + [2(c_2 - \sum_{k=i}^{r-1} l_k) - (2i-1)] + il_{i-1}. \blacksquare$$

**4.3. Generic bundles in the moduli space $\mathcal{M}_L(r; tf, c_2)$**

Our purpose is to determine the structure of a generic bundle in $\mathcal{M}_L(r; tf, c_2)$.

**Lemma 4.15.** *Assume $\mathcal{M}_L(r; tf, c_2)$ is nonempty where $-r < t \leq 0$ and $L$ satisfies (3.3). Then for a generic bundle $V$ in $\mathcal{M}_L(r; tf, c_2)$, there are $(r-1)$ exact sequences:*

$$0 \to \mathcal{O}_X(-n_i f) \to V_i \to V_{i-1} \to 0 \tag{4.16}$$

*for $r \geq i \geq 2$ with the following properties:*

(i) $V_r = V$, $V_i$ is a rank-$i$ bundle for $i = r-1, \ldots, 2$, and

$$V_1 = \mathcal{O}_X((t + \sum_{i=2}^{r} n_i)f) \otimes I_{Z_1};$$

(ii) $\ell(Z_1) = c_2$, and $Z_1$ is supported in $c_2$ distinct fibers;

(iii) $n_r = [\frac{c_2 - t}{r}]$, and $n_i = [\frac{(c_2 - t) - \sum_{k=i+1}^{r} n_i}{i}]$ for $i = r-1, \ldots, 2$.

*Proof.* Note that $\delta_1^{**} = \#(\text{moduli of } V_1^{**}) - \dim \text{Aut}(V_1^{**}) = -1$. By Proposition 4.14,

$$\delta_r^{**} \leq \delta_1^{**} + \sum_{i=2}^{r} [2c_2 - 2\sum_{k=i}^{r-1} l_k - (2i-1) + i l_{i-1}]$$

$$= -1 + [2(r-1)c_2 + (1 - r^2) + \sum_{i=1}^{r-1} (3-i) l_i].$$

Since $\delta_r^{**} = \#(\text{moduli of } V) - 1$ and $\sum_{i=1}^{r-1} l_i = c_2$, we have

$$\#(\text{moduli of } V) \leq 2rc_2 + (1 - r^2) + \sum_{i=1}^{r-1} (1-i) l_i \leq 2rc_2 + (1 - r^2). \tag{4.17}$$

By Lemma 3.6, since $\mathcal{M}_L(r; tf, c_2)$ is nonempty, we always have

$$\#(\text{moduli of } V) = 2rc_2 + (1 - r^2);$$

thus, in particular, all the inequalities in (4.17), (4.13) and Proposition 4.6 (iii) become equalities. Hence, for a generic bundle $V$ in $\mathcal{M}_L(r; tf, c_2)$, we conclude that:

(a) since (4.17) is an equality, $l_2 = \ldots = l_{r-1} = 0$; so $l_1 = c_2$. It follows that $V_2, \ldots, V_{r-1}$ are bundles, and (4.16) comes from (4.7). Since $V_1$ is of rank-1,

$$V_1 = \mathcal{O}_X((t + \sum_{i=2}^{r} n_i)f) \otimes I_{Z_1}$$

for some 0-cycle $Z_1$ on $X$. Thus, $Q_1 = \mathcal{O}_{Z_1}$, and $\ell(Z_1) = \ell_1 = c_2$. This proves (i).

(b) since (4.13) is an equality and $Q_1 = \mathcal{O}_{Z_1}$,

$$\#(\text{moduli of } Z_1) = \#(\text{moduli of } Q_1) = 2\ell_1 = 2c_2.$$

Thus, for a generic bundle $V$, $Z_1$ is reduced and supported in $c_2$ distinct fibers. This proves (ii).

(c) since Proposition 4.6 (iii) is an equality, for $i = 2, \ldots, r$, we have

$$i \cdot n_i + (i - j_i) = c_2(V_i^{**}) - t - \sum_{k=i+1}^{r} n_k = c_2 - t - \sum_{k=i+1}^{r} n_k;$$

note that $0 \leq (i - j_i) < i$; thus, $n_r = [\frac{c_2 - t}{r}]$, and

$$n_i = [\frac{(c_2 - t) - \sum_{k=i+1}^{r} n_i}{i}]$$

for $i = r - 1, \ldots, 2$. This proves (iii) and completes the proof. ∎

**Proposition 4.18.** *Assume that $\mathcal{M}_L(r; tf, c_2)$ is nonempty where $-r < t \leq 0$ and $L$ satisfies (3.3). Then, a generic bundle $V$ in $\mathcal{M}_L(r; tf, c_2)$ sits in an exact sequence:*

$$0 \to \bigoplus_{i=1}^{r} \mathcal{O}_X(-n_i f) \to V \to \bigoplus_{i=1}^{c_2} (\tau_i)_* \mathcal{O}_{f_i}(-1) \to 0 \tag{4.19}$$

*where the integer $n_i$ is defined by induction as follows:*

$$n_i = [\frac{(c_2 - t) - \sum_{k=i+1}^{r} n_i}{i}] \text{ for } i < r \text{ with } n_r = [\frac{c_2 - t}{r}], \tag{4.20}$$

*and $\{f_1, \ldots, f_{c_2}\}$ are distinct fibers with $\tau_i$ being the natural embedding $f_i \hookrightarrow X$.*

*Proof.* First of all, we notice that if $(c_2 - t) = ar + \epsilon$ with $0 \leq \epsilon < r$, then

$$n_i = \begin{cases} a & \text{if } i = \epsilon + 1, \ldots, r \\ a + 1 & \text{if } i = 1, \ldots, \epsilon. \end{cases} \tag{4.21}$$

In particular, $n_i \leq n_j$ if $i > j$. By Lemma 4.15, for a generic bundle $V$ in $\mathcal{M}_L(r; tf, c_2)$, we have $(r - 1)$ exact sequences (4.16). Consider the first two exact sequences:

$$0 \to \mathcal{O}_X(-n_r f) \to V \xrightarrow{p_{r-1}} V_{r-1} \to 0$$
$$0 \to \mathcal{O}_X(-n_{r-1} f) \to V_{r-1} \to V_{r-2} \to 0.$$

Then, the subsheaf $p_{r-1}^{-1}(\mathcal{O}_X(-n_{r-1}f))$ of $V$ sits in an exact sequence:

$$0 \to \mathcal{O}_X(-n_r f) \to p_{r-1}^{-1}(\mathcal{O}_X(-n_{r-1}f)) \to \mathcal{O}_X(-n_{r-1}f) \to 0.$$

Since $n_r \leq n_{r-1}$, $\mathrm{Ext}^1(\mathcal{O}_X(-n_{r-1}f), \mathcal{O}_X(-n_r f)) = 0$; thus,

$$p_{r-1}^{-1}(\mathcal{O}_X(-n_{r-1}f)) = \bigoplus_{i=r-1}^{r} \mathcal{O}_X(-n_i f).$$

We check that $V/\bigoplus_{i=r-1}^{r} \mathcal{O}_X(-n_i f) = V_{i-1}/\mathcal{O}_X(-n_{r-1}f) = V_{i-2}$. Thus, $V$ sits in

$$0 \to \bigoplus_{i=r-1}^{r} \mathcal{O}_X(-n_i f) \to V \to V_{r-2} \to 0.$$

By induction and the fact that $\mathrm{Hom}(\mathcal{O}_X(-n_1 f), V_1) \cong H^0(X, \mathcal{O}_X(c_2 f) \otimes I_{Z_1}) \neq 0$, we conclude that $V$ sits in an exact sequence:

$$0 \to \bigoplus_{i=1}^{r} \mathcal{O}_X(-n_i f) \to V \to V_1/\mathcal{O}_X(-n_1 f) \to 0.$$

Now, the exact sequence (4.19) follows from the observation that

$$V_1/\mathcal{O}_X(-n_1 f) = I_{Z_1}/\mathcal{O}_X(-c_2 f) = \bigoplus_{i=1}^{c_2} (\tau_i)_* \mathcal{O}_{f_i}(-1)$$

where $f_1, \ldots, f_{c_2}$ are the $c_2$ distinct fibers supporting the 0-cycle $Z_1$. ∎

Remark 4.22. (i) By Theorem 3.5, for any stable bundle $V$ in $\mathcal{M}_L(r; tf, c_2)$, $\pi^*(\pi_* V)$ is a locally free rank-$r$ subsheaf of $V$ with the quotient $Q$ being supported on the fibers of the ruling $\pi$ over which the restriction of $V$ is non-trivial. Another possible approach to prove Proposition 4.18 is to study the exact sequence

$$0 \to \pi^*(\pi_* V) \to V \to Q \to 0$$

and to estimate the number of moduli of these $V$'s in terms of the data of $Q$ and the rank-$r$ bundle $\pi_* V$ on $\mathbf{P}^1$. In fact, this approach has been used very successfully by Friedman [8] to study stable rank-2 bundles on an arbitrary ruled surface. However, for $r > 2$, the difficulty of this approach lies in the observation that the deformation of $Q$ is quite complicated.

(ii) From the exact sequence (4.19), we conclude that

$$\pi^*(\pi_*V) = \bigoplus_{i=1}^{r} \mathcal{O}_X(-n_i f)$$

for a generic bundle $V$ in the moduli space $\mathcal{M}_L(r; tf, c_2)$.

## 5. The moduli space $\mathcal{M}_L(r; tf, c_2)$ on a rational ruled surface

In this section, based on the results from the previous section, we determine the birational structure of the moduli space $\mathcal{M}_L(r; tf, c_2)$ on a rational ruled surface where $L$ satisfies (3.3) and $-r < t \leq 0$. First of all, we introduce the following notations.

*Notation 5.1.* (i) Let $n_i$, $f_i$ and $\tau_i$ be as in Proposition 4.18. Put

$$W_0 = \bigoplus_{i=1}^{r} \mathcal{O}_{\mathbf{P}^1}(-n_i), \quad W = \pi^*(W_0) = \bigoplus_{i=1}^{r} \mathcal{O}_X(-n_i f), \quad \text{and} \quad Q = \bigoplus_{i=1}^{c_2} (\tau_i)_* \mathcal{O}_{f_i}(-1);$$

(ii) Let $\mathcal{M}$ be the Zariski open and dense subset in $\mathcal{M}_L(r; tf, c_2)$ parametrizing all bundles sitting in exact sequences of the form (4.19);

(iii) Let $\Phi : \mathcal{M} \to U$ be the morphism defined by

$$\Phi(V) = \sum_{i=1}^{c_2} \pi(f_i)$$

where $U$ is a Zariski open and dense subset in $\mathrm{Sym}^{c_2}(\mathbf{P}^1) \cong \mathbf{P}^{c_2}$.

Next, we want to determine the fiber $\Phi^{-1}(u)$ for $u \in U$. We start with a lemma.

**Lemma 5.2.** (i) $\mathrm{Hom}(W, V) \cong \mathrm{End}(W)$;
(ii) $\dim \mathrm{Aut}(W) = r^2$ and $\dim \mathrm{Aut}(Q) = c_2$;
(iii) $\dim \mathrm{Ext}^1(Q, W) = 2rc_2$.

*Proof.* (ii) and (iii) follow from (4.21) and the definitions of $W$ and $Q$. In the following, we prove (i). Since $W = \pi^* W_0$, $\mathrm{End}(W) \cong \mathrm{End}(W_0)$. Since $\pi_* Q$ is torsion and

$$H^0(\mathbf{P}^1, \pi_*Q) = H^0(X, Q) = 0,$$

$\pi_*Q$ must be zero. Applying $\pi_*$ to (4.19), we have $\pi_*V \cong \pi_*W = W_0$. Thus,

$$\begin{aligned}\mathrm{Hom}(W,V) &\cong H^0(X, V \otimes W^*) = H^0(\mathbf{P}^1, \pi_*(V \otimes \pi^*(W_0^*))) \\ &\cong H^0(\mathbf{P}^1, W_0 \otimes W_0^*) \cong \mathrm{End}(W_0) \\ &\cong \mathrm{End}(W). \blacksquare\end{aligned}$$

**Proposition 5.3.** *Let $u \in U$. Then, the fiber $\Phi^{-1}(u)$ is birational to $\mathrm{Ext}^1(Q,W)$ modulo the $(c_2 + r^2 - 1)$-dimensional group actions from $\mathrm{Aut}(W)/\mathbf{C}^*$ and $\mathrm{Aut}(Q)$.*

*Proof.* By Lemma 5.2 (i), $\mathrm{Hom}(W,V) \cong \mathrm{End}(W)$. From the proof of Lemma 4.15, we see that generic extensions in $\mathrm{Ext}^1(Q,W)$ must correspond to bundles in the Zariski open and dense subset $\mathcal{M}$. It follows that $\Phi^{-1}(u)$ is birational to $\mathrm{Ext}^1(Q,W)$ modulo the group actions from $\mathrm{Aut}(W)/\mathbf{C}^*$ and $\mathrm{Aut}(Q)$. By Lemma 5.2 (ii),

$$\dim \mathrm{Aut}(W) = r^2 \quad \text{and} \quad \dim \mathrm{Aut}(Q) = c_2.$$

Therefore, the group actions are $(c_2 + r^2 - 1)$-dimensional. $\blacksquare$

Now, we prove the second main result in this paper.

**Theorem 5.4.** *Assume that the moduli space $\mathcal{M}_L(r;tf,c_2)$ is nonempty where $r \geq 2$, $-r < t \leq 0$ and the ample divisor $L$ satisfies the condition (3.3). Then,*

(i) *$\mathcal{M}_L(r;tf,c_2)$ is irreducible and unirational;*

(ii) *a generic bundle $V$ in $\mathcal{M}_L(r;tf,c_2)$ sits in an exact sequence*

$$0 \to \bigoplus_{i=1}^r \mathcal{O}_X(-n_i f) \to V \to \bigoplus_{i=1}^{c_2} (\tau_i)_* \mathcal{O}_{f_i}(-1) \to 0 \tag{5.5}$$

*where the integer $n_i$ is defined by induction as follows:*

$$n_i = \left[\frac{(c_2 - t) - \sum_{k=i+1}^r n_i}{i}\right] \text{ for } i < r \text{ with } n_r = \left[\frac{c_2 - t}{r}\right], \tag{5.6}$$

*and $\{f_1, \ldots, f_{c_2}\}$ are distinct fibers with $\tau_i$ being the natural embedding $f_i \hookrightarrow X$;*

(iii) *$(c_2 - t) \geq r$.*

*Proof.* (i) By Lemma 5.2 (iii), the extension group $\mathrm{Ext}^1(Q,W)$ has dimension $2rc_2$. By Proposition 4.24, we have a rational map $\Phi$ from the moduli space $\mathcal{M}_L(r;tf,c_2)$ to $\mathbf{P}^{c_2}$ such that a generic fiber $\Phi(u)$ is birational to

$$[\mathrm{Aut}(W)/\mathbf{C}^*] \backslash \mathbf{C}^{\oplus 2rc_2} / \mathrm{Aut}(Q).$$

Therefore, $\mathcal{M}_L(r; tf, c_2)$ is irreducible and unirational.

(ii) This is the same as Proposition 4.18.

(iii) Since $\mathcal{O}_X(-n_r f) \hookrightarrow V$ and $V$ is $L$-stable, $-n_r f \cdot L < tf \cdot L/r \leq 0$; thus, $n_r \geq 1$. Since $n_r = [(c_2 - t)/r] \leq (c_2 - t)/r$, we get $(c_2 - t) \geq r$. ∎

*Remark 5.7.* In the Theorem 1.9 of [2], Artamkin showed that if $c_2 \geq r \geq 2$, then $\mathcal{M}_L(r; 0, c_2)$ is nonempty and irreducible. Therefore, by Theorem 5.4 (iii), we conclude that $\mathcal{M}_L(r; 0, c_2)$ is nonempty if and only if $c_2 \geq r$.

Department of Mathematics, Hong Kong University of Science and Technology, Clear Water Bay, Kowloon, Hong Kong; E-mail address: mawpli@uxmail.ust.hk

Department of Mathematics, Oklahoma State University, Stillwater, OK 74078, U.S.A. E-mail address: zq@hardy.math.okstate.edu.